\def \bgamma{\bar{\gamma}}
\def \bomega{\bar{\omega}}
\def \tK{E}
\def \bet{\bar{\eta}}
\def \bxi{\bar{\xi}}
\def \balpha{\bar{\alpha}}
\def \bbeta{\bar{\beta}}
\def \bx{\bar{x}}
\def \o{\omega}
\def \bU{\bar{U}}
\def \bs{\bar{\s}}
\def \sumn{\sum_{n=1}^g}
\def \summ{\sum_{m=1}^g}
\def \be{\bar{e}}
\def \bK{\bar{K}}
\def \X{{\bf X}}

\def \bA{\bar{A}}
\def \bl{\bar{\lambda}}
\def \ba{\bar{a}}

\def \bpsi{\bar{\psi}}
\def \O{{\cal O}}

\def \bs{\bigskip}
\def \st{\stackrel}

\def \d{\delta}

\def \s{\sigma}
\def \a{\alpha}
\def \b{\beta}

\def \e{\epsilon}
\def \p{^\prime}

\def \Proof{\medskip\noindent{\bf \it Proof.\quad}}

\def \Ga{{\bf G}_a}
\def \ra{\rightarrow}

\documentstyle [12pt]{article}
\normalbaselines
\vspace{1.1in}
\title{Reduction of the Manin map modulo $p$}
\author{Alexandru Buium and Jos\'e Felipe Voloch}
\date{\  }
\begin{document}
\maketitle

\bs

  For an abelian variety $A$ over a function field $K$ of characteristic
zero equipped with a derivation $\d:K \ra K$ Manin defined in [Man1], [Man2]
a remarkable additive map $A(K) \ra V$, where $V$ is a  vector space over $K$,
which plays an important role in
diophantine geometry over function fields. (Cf. [Co] for a ``modern" exposition
of Manin's work. Cf. also [B1], [B2] for a different way of introducing
this map.) In the ``generic case" this map is a ``second order non linear
differential operator". In [V1] the second author defined an analogue of this
map in the case of elliptic curves over function fields of characteristic $p$.
This analogous map turned out to be of order one.
 Then the following  results were proved in [V1] for ordinary elliptic curves:

\bs

\noindent (*) {\it  The ``reduction mod} $p$" {\it of the Manin map in
characteristic zero is the ``derivative" of the Manin
map in characteristic} $p$.

\bs

\noindent (**) {\it The kernel of the Manin map in characteristic} $p$ {\it is
the group
of points divisible by} $p$.

\bs

 In the higher dimensional case one also
has a ``characteristic $p$ analogue of the Manin map", which is implicit in
[BV].
(Cf. also [H] for a different approach.) The aim of this paper is to prove
(*) and (**) in the higher dimensional case.

The method used in [V1] to prove (*) for elliptic curves was based on a
computation with Tate curves hence  seems to be hard to extend to the higher
dimensional case; instead
we use the approach in [B1], [B2].
Note also that since (*) is a statement about algebraic groups (in
characteristic $p >0$), it makes sense
to consider the statement which (*) implies at the level of Lie algebras;
it turns out that the corresponding Lie algebra statement is (a dual form
of) a basic result of Igusa-Manin-Katz [Man3], [K] which says, roughly
speaking, that the Hasse-Witt matrix satisfies the Picard-Fuchs equation.
But of course there is no way back, in characteristic $p > 0$, from Lie
algebras to algebraic groups; so there is no way back from the Igusa-Manin-Katz
result
to our result.

Note on the other hand that (**) is an analogue in characteristic $p$ of the
Manin-Chai ``Theorem of the kernel" [Man1] [Man2], [Ch];
we shall prove (**) by relating  the Manin map in characteristic
$p$ to the $p$-descent map.

\bs

 To state our main result let us recall first some general definitions  from
[B1] which lead to a general concise  definition of the ``Manin maps" in
arbitrary characteristic.
 Let $A$ be a scheme over a scheme $S$ and assume we are given a derivation
$\d$ on the structure sheaf of $S$. Then one can form a projective system of
$S-$schemes $(A^n)$ for $n \geq -1$ with affine transition maps $\pi_n$
and derivations $\d=\d_n:\O_{A^n} \ra \pi_{n*}\O_{A^{n+1}}$ as follows.
Set $A^{-1}=S$, $A^0=A$;  let $\d_{-1}$ be induced by $\d$ and define
inductively
$A^{n+1}=Spec\ (S(\Omega_{A^n})/I_n)$ where $I_n$ is the ideal generated by
sections of the form $d(\pi_{n-1}^*f) - \d_{n-1}f,\ f \in \O_{A^{n-1}}$,
while $\d_n$ is induced by the Kahler differential $d$ in the obvious way.
The schemes $A^n$ are called the
{\it schemes of} $n-${\it jets of} $A/S$ {\it along the direction} $\d$,
or simply the {\it canonical prolongations of} $A/S$. This construction
commutes, in the obvious sense, with ``horizontal" base change $(S\p, \d\p) \ra
(S,\d)$ (here ``horizontal"
means ``compatible with the derivations") and has the following universality
property: for any $A^{n}-$scheme $Z$ and any derivation $\partial$
from $\O_{A^n}$ to (the direct image of) $\O_Z$ prolonging the derivation
$\d_{n-1}$ there exists a unique $A^n-$scheme morphism $Z \ra A^{n+1}$
such that $\partial$ is induced by $\d_{n}$. We refer to [B2] Part I for
details.

Assume now in addition that $A/S$ is a smooth group scheme. Then $(A^n)$ will
be a projective system of smooth group schemes. Denote by $\X^n(A)$
the set of all $S-$group scheme homomorphisms from $A^n$ to the additive group
$\Ga$ over $S$. This set has actually a structure of $\O(S)-$submodule of the
ring $\O(A^n)$. Moreover the maps $\O(A^n) \ra \O(A^{n+1})$ induced by the
$\pi_n$'s are injective and will be viewed as inclusions; so we get induced
derivations $\d:\O(A^n) \ra \O(A^{n+1})$ which induce maps $\d:\X^n(A) \ra
\X^{n+1}(A)$. The space $\X^n(A)$ was called in [B2] the {\it space of}
$\d-${\it polynomial characters of} $A$ {\it of order} $\leq n$. Note that
each element $\psi \in \X^n(A)$ defines a
homomorphism $\hat{\psi}:A(S) \ra \O(S)$ by the formula
\[\hat{\psi}(P)=\psi_S(P^n),\ \ \ P \in A(S)\]
where $\psi_S:A^n(S) \ra \Ga (S) =\O(S)$ is the map induced by $\psi$
and $P^n \in A^n(S)$ stands for the ``canonical lifting" of $P$ (induced by the
universality property).
The components of the ``classical Manin map" in characteristic zero [Man1]
[Man2] as well as the
Manin map in characteristic $p$ in [V1] are all of the form $\hat{\psi}$ above;
so studying ``Manin maps" is the same as studying the spaces $\X^n(A)$.

\bs

 Throughout the paper we shall consider the following situation. We start with
a discrete valuation ring $R$
 and we denote by
$K$ and $\bK$ the quotient field and the residue field of $R$ respectively.
We assume $K$ has characteristic zero and $\bK$ is a function field of
one variable over a perfect field of characteristic $p >0$.
As a general rule the upper bar will denote the reduction modulo the maximal
ideal $m_R$ of $R$; in particular for any element $a \in R$ we denote its image
in $\bK$  by $\ba$.
We assume we are given a derivation $\d:R \ra R$ such that $\d (m_R) \subset
m_R$ (this is automatic if $p$ is tamely ramified in $R$) and such that $\d (R)
\not\subset m_R$. Then this derivation will induce a non zero derivation (still
denoted by) $\d$ on $\bK$.
(A typical example
of this situation is: $R={\bf Z}[t]_{(p)}$, $\d=\partial/\partial t$,
$K={\bf Q}(t)$, $\bK={\bf F}_p(t)$.)

Next we consider an abelian scheme $A/R$ of relative dimension $g$; let $A_K/K$
and $\bA / \bK$ be the generic and special fibres
respectively and consider the corresponding spaces of $\d-$polynomial
characters
$\X^n(A), \X^n(A_K), \X^n(\bA)$.
We will make the following two assumptions:

\bs

(i) $\bA$ has $p-$rank $g$ (one also says that $\bA$ is {\it ordinary}).
Recall that this means that the rank of the Frobenius endomorphism of
$H^1(\bA,\O)$
equals $g$ (recall that the rank of a $p-$linear map is the dimension
of the linear span of its image).

\bs

(ii) $\bA$ has $\d-$rank $g$ (cf. [B2], Part I). Recall that this means that
the $\bK-$linear map $\rho(\d) \cup: H^0(\bA, \Omega^1) \ra H^1(\bA,\O)$
induced by cup product with the Kodaira Spencer class $\rho(\d) \in H^1(\bA,
T)$ has rank $g$ (we denoted by $\rho:Der\ \bK \ra H^1(\bA, T)$ the
Kodaira-Spencer map).

\bs

Of course (ii) implies that $A_K$ itself has $\d-$rank $g$ (in the analogous
sense). So by [B2], Part I,  $\X^1(A_K)=0$ and that $\X^2(A_K)$ has dimension
$g$ over $K$. The ``classical Manin map in characteristic zero" in this case is
the map
$(\hat{\psi}_1,...,\hat{\psi}_g):A(K) \ra K^g$ associated to a basis
$\psi_1,...,\psi_g$ of $\X^2(A_K)$.

 We have a natural inclusion
$\X^2(A) \subset \X^2(A_K)$; elements of $\X^2(A_K)$ which lie in $\X^2(A)$
will be called {\it integral}. We also have  a reduction modulo $m_R$ map
\[\X^2(A) \ra \X^2(\bA),\ \ \ \psi \mapsto \bpsi\]
so the integral elements $\psi$ of $\X^2(A_K)$ may be reduced modulo $m_R$ to
get elements $\bpsi \in \X^2(\bA)$. Finally recall that we have  have a map
induced by derivation
\[\X^1(\bA) \ra \X^2(\bA),\ \ \  \phi \mapsto \d \phi\]
Our main result is the following:

\bs

\proclaim Theorem. There exists an integral basis $\psi_1,...,\psi_g$
of $\ \X^2(A_K)$
 and elements $\phi_1,...,\phi_g \in \X^1(\bA)$
 with the following properties:\\
\ \\
1) $\bpsi_n=\d \phi_n$ for $1 \leq n \leq g$.\\
\ \\
2) $\{ P \in \bA(\bK) | \hat{\phi_1}(P)=\dots=\hat{\phi_g}(P)=0\}=
p\bA(\bK)$.

\bs

Actually, as we shall recall below, there is a natural increasing ``filtration
by degrees"
$F^d \X^1(\bA), d \geq 0$ on $\X^1(\bA)$ (cf. [B2], Part II); then we shall
prove that
$F^{p-1} \X^1(\bA)=0$ and that  $\phi_1,...,\phi_g$ in the Theorem may be taken
to be a basis of $F^p \X^1(\bA)$.

\bs

In what follows we devote ourselves to the proof of the Theorem. In the end of
the paper we will show how to deduce from our Theorem the (dual form of the)
Igusa-Manin-Katz
theorem by passing to Lie algebras.

\bs

 We need a ``cocycle description" of the various objects involved; cf. [B3].
Let $U_i$ be an affine covering of $A$ and let $\theta_i$ be derivations
of $\O(U_i)$ which lift the derivation $\d$ of $R$. Fix a basis
$\o_1,...,\o_g$ of the $R-$module  $H^0(A,\Omega^1)$ and let
$v_1,...,v_g$ be the dual basis of  $L(A)=$Lie algebra of $A/R$. We view
elements
of $L(A)$ as derivations of $\O_A$.  Then we may write
\[(1)\ \ \ \theta_j - \theta_i = \sumn a_{ijn}v_n \]
with $a_{ijn} \in \O(U_i \cap U_j)$. Then the $a_{ijn}$ are cocycles; let
$e_n \in H^1(A,\O)$ be the classes of these cocycles. Since the reduction
modulo $m_R$
 of the cocycle $\theta_j - \theta_i$ represents the Kodaira Spencer class
$\rho(\d) \in H^1(\bA,T)$ it follows that the images $\be_1,...,\be_g \in
H^1(\bA,\O)$ of the $e_i$'s are a basis of $H^1(\bA,\O)$
(image of $\bomega_1,...,\bomega_g \in H^0(\bA, \Omega^1)$ via the map
$\rho(\d): H^0(\bA, \Omega^1) \ra H^1(\bA, \O)$)
hence $e_1,...,e_g$
form a basis of $H^1(A,\O)$.

Now our choice of a basis for the tangent bundle of $U_i$ and of the lifting
$\theta_i$ provides an $U_i-$isomorphism
\[\s_i:U_i^1 \ra U_i \times Spec\ R[x_1,...,x_g] \]
such that the derivation $\d:\O(U_i) \ra \O(U_i^1)$ corresponds to the
derivation
\[ (2)\ \ \ \ \d_i:\O(U_i) \ra \O(U_i)[x_1,...,x_g],\ \ \d_i=\theta_i + \sumn
x_nv_n\]
It follows from the general theory in [B2] Part I that $A^1$ is the universal
vectorial extension of $A$ by $\Ga^g$ (actually in [B2] one assumes a ground
field rather than a ground discrete valuation ring, but all arguments go
through). Then
the isomorphisms $\s_i$ are $\Ga^g-$equivariant where $\Ga^g$ acts on the
affine space $Spec\ R[x_1,...,x_g]$ by translations on the affine coordinates.
Setting $x_{in}:=\s_i^*x_n$ we get from (1) and (2) that
\[(3)\ \ \ \ \ \ \ a_{ijn}=x_{in}-x_{jn}\]
in the ring $\O(U_{ij}^1)$, $U_{ij}:=U_i \cap U_j$. Applying $\d$ we get
\[(4)\ \ \ \ \ \ \ \d  a_{ijn}=\d x_{in}- \d x_{jn}\]
in the ring $\O(U_{ij}^2)$.

Reducing  the isomorphisms $\s_i$ modulo $m_R$ we get isomorphisms
\[ \bar{\s}_i:\bU_i^1 \ra \bU_i \times Spec\ \bK[\bx_1,...,\bx_g] \]
The pull backs of $\bx_1,...,\bx_g$ via $\bar{\s}_i$ restricted to
$L(\bA)=Ker(\bA^1 \ra \bA)$ will be affine maps whose associated linear maps
are $\bomega_1,...,\bomega_g$.
Now the rings $\O(\bU_i)[\bx_1,...,\bx_g]$ have a natural filtration by degrees
(the $d-$piece of the filtration being the space of all polynomials
of degree $\leq d$). These filtrations induce via $ \bar{\s}_i$ filtrations
on $\O(\bU_i^1)$ (the $d-$th piece of this filtration is the space
of all $\bK-$linear combinations of $e-$fold products of elements
$\bx_{i1},...,\bx_{ig}$ with $e \leq d$.) Clearly these filtrations glue
together
to give a filtration $F^d \O(\bA^1)$ on $\O(\bA^1)$ and hence an induced
filtration $F^d \X^1(\bA)$ on $\X^1(\bA)$. We could introduce similar
filtrations on $\X^1(A), \X^1(A_K)$ but we won't need them in what follows.
Here are two basic facts about the above introduced filtration.

\bs

\proclaim Lemma 1. $F^{p-1}\X^1(\bA)=0$.

\bs

\Proof Let $\phi \in   F^{p-1}\X^1(\bA)$. Then the restrictions of $\phi$ to
$\bU_i^1$ have the form $\bar{\s}_i^*H_i$ where $H_i=H_i(\bx_1,...,\bx_g)$ are
polynomials of degree $\leq p-1$ with coefficients in $\O(\bU_i)$. By the
$\Ga^g-$equivariance of $\s_i$, $H_i$ must have the form
\[H_i=h_i+ \sum_{m \geq 0} \sumn h_{inm} \bx_n^{p^m}\]
where $h_i, h_{inm} \in \O(\bU_i)$.
Since we are looking at the $F^{p-1}-$piece of the filtration we must
have $h_{inm}=0$ for $m \geq 1$. The conditions $(\bar{\s}_j^*)^{-1}
\bar{\s}_i^*H_i = H_j$ give
\[(5)\ \ \ \ h_i + \sumn h_{in0}(\bx_n + \ba_{ijn})=h_j + \sumn h_{jn0} \bx_n\]
This implies that  $(h_{in0})_i$ glue together to give
an element $h_n \in \O(\bA)=\bK$. Then (5) further implies that
\[\sumn h_n \ba_{ijn} = h_j - h_i\]
Passing to cohomology classes one gets $\sumn h_n \be_n =0$ in $H^1(\bA, \O)$
so we get $h_n=0$ for all $n$. So $H_i=h_i$ for all $i$. Again the $h_i$ glue
together to give an element in $\bK$ which must be $0$ hence $\phi=0$ and we
are done.

\bs

\proclaim Lemma 2. dim\  $F^{p}\X^1(\bA)=g$.

\bs

\Proof First we ``explicitely" construct $g$ linearly independent elements
$\phi_1,...,\phi_g$ in $F^{p}\X^1(\bA)$; this construction will be used later.
By the $p-$rank condition we may write
\[\be_n^{(p)} = \summ \bl_{nm} \be_m,\ \ \ 1 \leq n \leq g\]
$\lambda_{nm} \in R,\  det(\bl_{nm}) \neq 0$, where the upper $(p)$ means
``image under the Frobenius". The matrix $\bl=(\bl_{nm})$ is classically called
the {\it Hasse-Witt matrix}  (corresponding to the basis
$\bar{e}_1,...,\bar{e}_g$). So we may write
\[(6)\ \ \ \ \ba_{ijn}^p - \summ \bl_{nm} \ba_{ijm} = \ba_{in} - \ba_{jn},\ \ 1
\leq n \leq g\]
where $a_{in} \in \O(U_i)$. Define
\[(7)\ \ \ \ \phi_{in}:= \bx_{in}^p - \summ \bl_{nm} \bx_{im} - \ba_{in} \in
\O(\bU_i^1)\]
Then due to (6) the $(\phi_{in})_i$ glue together to give an element
$\phi_n \in \O(\bA^1)$. Subtracting for each $n$ an element  $\mu_n \in R$ from
all the $a_{in}$'s we may assume $\phi_n(0)=0$ for all $n$. Clearly
$\phi_1,...,\phi_g$ are $\bK-$linearly independent. Let us check that they are
elements of $F^{p}\X^1(\bA)$; we only have to check that they are additive
characters on $\bA^1$. Let $V$ be the kernel of $\bA^1 \ra \bA$. We have
\[(8)\ \ \ \phi_n(u+v)=\phi_n(u)+\phi_n(v),\ \ \ u \in \bA^1, v \in V\]
because $\bar{\s}_i$ is $\Ga^g-$equivariant and transforms the right hand side
of (7) into an ``affine polynomial" (i.e. an additive polynomial plus
a term of degree zero). But property (8) immediately implies that $\phi_n$
is additive; indeed for any fixed $u \in \bA^1$ the regular function on $\bA^1$
defined by $u \mapsto \phi_n(u+v) - \phi_n(u) - \phi_n(v)$ vanishes at $0$ and
is constant on the fibres of $\bA^1 \ra \bA$. Since $\O(\bA)=\bK$ the above
function is $0$, hence $\phi_n$ is additive.

Now, exactly as in Lemma 1, any element $\phi \in F^{p}\X^1(\bA)$ is
represented by polynomials $H_i$ of the form
\[H_i=h_i + \sumn h_{in0} \bx_n + \sumn h_{in1} \bx_n^p\]
where $h_i, h_{in0}, h_{in1} \in \O(\bU_i)$. As in Lemma 1 one gets that
$(h_{in1})_i$ glue together to give an element $h_{n1} \in \bK$. Hence
\[\phi - \sumn h_{n1} \phi_n \in F^1\X^1(\bA)\]
so by Lemma 1 $\phi - \sumn h_{n1} \phi_n=0$ and we conclude that
$\phi_1,...,\phi_g$ generate $ F^{p}\X^1(\bA)$ which concludes the proof of
Lemma 2.

\bs

Let us make the remark (to be used later) that the two Lemmas above do not
depend on the fact that our characteristic $p$ situation lifts to
characteristic zero. In particular they hold if we replace $\bK$ by a finite
separable extension $\tK$ of it. (We do not need to assume that the derivation
on $\tK$
lifts to a situation in characteristic zero.)

\bs

                                                                       Let us
come back to the proof of our Theorem.
One of the key steps   is the following remark: we know from [B3], Lemma
(2.10), p. 73 that the derivation $\d$ of $\O_A$ into the direct image of
$\O_{A^1}$ lifts to a derivation $\tilde{\d}$ of the whole of $\O_{A^1}$ ,
which, via the isomorphisms
$\s_i$ corresponds to derivations
\[ \tilde{\d}_i:\O(U_i)[x_1,...,x_g]
 \ra \O(U_i)[x_1,...,x_g] \]
given by formulae of the form
\[\tilde{\d}_i=\theta_i + \sumn x_nv_n + \sumn L_{in}(x_1,...,x_g)
\frac{\partial}{\partial x_n} \]
where $L_{in}$ are polynomials of degree $\leq 1$ in $x_1,...,x_g$ with
coefficients
in $\O(U_i)$. This result was proved in [B3] for a ground field rather than for
 a ground valuation ring but the same arguments go through in our situation.
Note that the existence of the lifting $\tilde{\d}$ is actually a consequence
of the Grothendieck-Messing-Mazur theory [MM]; but the additional information
that $deg\ L_{in} \leq 1$ provided by [B3] will be crucial below !

By the universality property of canonical prolongations, $\tilde{\d}$ induces
a section $s:A^1 \ra A^2$ of the projection $\pi_2:A^2 \ra A^1$ such that
the corresponding map $s^*:\O(U_i^2) \ra \O(U_i^1)$ maps $\d x_{in}$
into
$L_{in}(x_{i1},...,x_{ig})$. Since the map $s^*:\O(U_{ij}^2) \ra \O(U_{ij}^1)$
is the identity on $\O(U_{ij}^1)$ applying this map to (4) we get
\[(9) \ \ \ \ \ \ \ \d  a_{ijn}= b_{in} - b_{jn},\ \ \
b_{in}=L_{in}(x_{i1},...,x_{ig}) \]

\bs

In what follows we shall construct an integral basis $\psi_1,...,\psi_g$ of
$\X^2(A)$. Let $(\lambda_{nm})$ be the matrix appearing in the proof of Lemma 2
and define
\[ (10)\ \ \ \ \psi_{in}:=\summ \lambda_{nm} (- \d x_{im} + b_{im}) \in
\O(U_i^2)\]
By (9) $(\psi_{in})_i$ glue together to give an element $\psi_n \in \O(A^2)$.
Subtracting, for each $n$, an element $\nu_n \in R$ from all of the $b_{in}$'s
we may assume that $\psi_n(0)=0$ for all $n$. Clearly $\psi_1,...,\psi_g$
are $K-$linearly independent. An argument similar to the one in the proof of
Lemma 2 shows that $\psi_n$ are additive characters (instead of $\O(\bA)=\bK$
one uses the fact that $\O(A_K^1)=K$).

\bs

To complete the proof of the first part of the
Theorem  we will check that, with $\phi_1,...,\phi_g$ as in the proof of Lemma
2 and with $\psi_1,...,\psi_g$ as in (10) above we have
$\bar{\psi}_n=\d \phi_n$ for all $n$. Now
\[\bar{\psi}_{in}-\d \phi_{in}= \sum \bl_{nm} \bar{b}_{im} + \sum (\d
\bl_{nm})\bx_{im} + \d \ba_{in}\]
By (9) we have
\[\bar{b}_{im}=\bar{L}_{im}(\bx_{i1},...,\bx_{ig})\]
while on the other hand by (2) we have
\[\d \ba_{in} = \theta_i \ba_{in} + \summ (v_m \ba_{in})x_{im}\]
so we see that
\[\bar{\psi}_{n}-\d \phi_{n} \in F^1\X^1(\bA)\]
By Lemma 1 we get $\bar{\psi}_{n}-\d \phi_{n}=0$ which completes the proof of
assertion 1).

\bs

To prove the second part of the theorem we will relate the
Manin map in characteristic $p$ to the $p$-descent
map, by the following construction, which generalizes that of [V1].
The isogeny of $\bA$ to itself  defined by ``multiplication by $p$" factors as
$
V \circ F$ where $F: \bA \to \bA^{(p)}$ is
the Frobenius and $V: \bA^{(p)} \to \bA$ is the Verschiebung. Since the latter
i
s etale, the points of $\ker V$ are rational over the separable closure $\tK$
of $\bK$. Clearly, $\{ x \in E|  \d x =0\}=E^p$.
 By
Cartier duality $\ker F$ is $\tK-$isomorphic to $\mu_p^g$, hence
$$H^1(\tK, \ker F) = H^1(\tK,\mu_p^g) = (\tK^*/(\tK^*)^p)^g \hookrightarrow
\tK^g,$$
where
$H^1$ stands for the first flat cohomology group of group schemes and
 the last map is induced by the logarithmic derivative $\tK^* \to \tK,
x \mapsto \delta x/x$.

Now, the coboundary map in flat cohomology $\bA^{(p)}(\tK) \to H^1(\tK, \ker F)
=(\tK^*/(\tK^*)^p)^g$ can be
given by $P \mapsto (f_1(P),\ldots,f_g(P))$, where the functions $f_1,
\ldots,f_g$ have divisors $pD_1, \ldots, pD_g$, such that $D_1, \ldots, D_g$
form a basis for the $p$-torsion of the Picard variety of $\bA^{(p)}$.
Indeed, given such a divisor $D_i$, we get a map from ${\bf Z}/p{\bf Z}$
to the $p$-torsion of the Picard variety of $\bA^{(p)}$ and by duality,
a map from $\ker F$ to $\mu_p$. Also $f_i \circ [p] = g_i^p$ for some
function $g_i$ and the map $H^1(\tK, \ker F) \to (\tK^*/(\tK^*)^p)$ is
given by associating the torsor $F^{-1}(P)$ of $\ker F$ to the torsor
$g_i(F^{-1}(P))$ of $\mu_p$ and, clearly, this torsor corresponds to $f_i(P)$.

Therefore the composite map $\beta: \bA^{(p)}(\tK) \to \tK^g$ is given by
$\delta$-polynomial characters of order 1 and degree 1. Also, by construction
the kernel of $\beta$ is $F(\bA(\tK))$.
Now, the proof of the proposition in section 4 of [V2] shows that the matrix
for
med by the
images of an ${\bf F}_p$-basis of $\ker V$ in $(\tK^*/(\tK^*)^p)^g$ is
essentially the matrix of the
Serre-Tate parameters of $\bA$ (modulo $p$-th powers)
at any place of $\bK$ of good, ordinary, reduction for $\bA$.
Also, a theorem of Katz shows that its image in $E^g$ is the matrix of the
cup-product with the Kodaira-Spencer class and therefore, is of maximal rank.
(This argument is given in detail in [V2]).
Let $\wp: \tK^g \to \tK^g$ be an additive polynomial map whose kernel
is the image of $\ker V$; then $\wp$ is clearly of degree $p$.
We now define $\mu: \bA(\tK) \to \tK^g$
by $\mu(P) = \wp(\beta(Q))$, where $Q \in \bA^{(p)}(E)$ is such that $V(Q) =
P$.
It is clear that the definition is independent of the choice of $Q$ and is
given
by $\delta$-polynomial characters of order 1 and degree $p$ on $\bA$.
By Lemma 2 (and the remark after it) the components of $\mu$ are defined by
$\tK-$linear combinations of $\phi_1,...,\phi_g$. Let us show
that the kernel of $\mu$ is $p\bA(\tK)$. Suppose $\mu(P) = 0$, then
after changing $Q$ by an element of $\ker V$, we get $\beta(Q) = 0$, so
$Q = F(R), R \in \bA(E)$, hence $P=V(F(R))=p R$.

Now the condition $Ker\ \mu=p \bA(\tK)$ implies that
\[\{ P \in \bA(\tK) | \hat{\phi_1}(P)=\dots=\hat{\phi_g}(P)=0\}=
p\bA(\tK) \]
We get that
\[\{ P \in \bA(\bK) | \hat{\phi_1}(P)=\dots=\hat{\phi_g}(P)=0\}=
p\bA(\tK) \cap \bA(\bK)  \]
We claim that
\[p\bA(\tK) \cap \bA(\bK) = p \bA(\bK) \]
which will close the proof. To check the claim note that the proposition in
[V2]
 says that $\bA(E)$ has no point of order $p$; so if $P=p Q \in \bA(\bK)$ for
so
me $Q \in \bA(E)$, then for all $\s \in Gal(E/\bK)$ we have $p(Q-Q^{\s})=0$
henc
e $Q-Q^{\s}=0$ so $Q \in \bA(\bK)$ and we are done.

\bs

{\it Remark.} Assume we are in the hypothesis of the Theorem.
Let
\[\nabla_{\d}:H^1_{DR}(A) \ra H^1_{DR}(A)\]
be the additive map  obtained by evaluating the Gauss-Manin connection at $\d$
and view $H^0(A, \Omega^1)$ as embedded into $H^1_{DR}(A)$.
Let $\o=(\o_1,...,\o_g)$ be an $R-$basis of $H^0(A, \Omega^1)$
and   write the Picard Fuchs equation:
\[(11) \ \ \ \ \ \ \nabla_{\d}^2 \o + \a \nabla_{\d} \o + \b \o =0 \]
where $\a, \b$ are $g \times g$ matrices with entries in $R$. Moreover
let $\be_1,...,\be_g$ be the image of $\bomega_1,...,\bomega_g \in
H^0(\bA,\Omega^1)$ via the isomorphism $\rho(\d): H^0(\bA,\Omega^1) \ra
H^1(\bA,\O)$
given by the cup product with the Kodaira Spencer class $\rho(\d) \in
H^1(\bA,T)$ and consider the Hasse-Witt matrix $\bl$ with respect to the basis
$\be$, in other words write $\be^{(p)}=\bl \be$. Then we claim that $\bl$
satisfies the following ``dual Picard-Fuchs equation":
\[(12)\ \ \ \ \ \ \d^2 \bl - (\d \bl) \bar{\a} + \bl (\bar{\b} - \d
\balpha)=0\]
The above claim is a dual version of the Igusa-Manin-Katz theorem referred to
in the beginning of the paper.

To check the claim above, let $\psi_1,...,\psi_g$ be as in the Theorem and let
us borrow the notations from the proof of the Theorem. The $\psi_i$'s define an
$R-$group scheme homomorphism $\psi: A^2 \ra \Ga^g$. We have induced maps
\[H^0(A,\Omega^1)=L(A)^{\circ} \subset L(A^2)^{\circ} \st{\psi^*}{\leftarrow}
L(\Ga^g)^{\circ}=R^g\]
where the upper $\circ$ means ``dual". Recall from [B2] Part III that one has
natural induced derivations
\[L(A)^{\circ} \st{\d}{\ra} L(A^1)^{\circ} \st{\d}{\ra} L(A^2)^{\circ}
\st{\d}{\ra} ...\]
and that there exists a basis $z=(z_1,...,z_g)$ of $R^g$
such that
\[\psi^*z=\d^2 \o + \a \d \o + \b \o\]
where $\a, \b$ are the matrices appearing in (11).
(Recall that this was deduced in [B2] Part III as  a consequence of the
Grothendieck-Messing-Mazur theory [MM], cf. also [B3], Chapter 2).
Moreover $\o, \d \o$ form a basis of $L(A^1)^{\circ}$ while $\o, \d \o, \d^2
\o$ form a basis of $L(A^2)^{\circ}$.

Now the Theorem implies that the map $\bar{\psi}^*:\bK^g \ra L(\bA^2)^{\circ}$
agrees with the composition
\[\bK^g \st{\phi^*}{\ra} L(\bA^1)^{\circ} \st{\d}{\ra} L(\bA^2)^{\circ}\]
on the standard basis $\e = ((1,0,...,0),(0,1,...,0),...)$ of $\bK^g$.
Write $\bar{z}=\bar{\xi} \e$ where $\xi$ is an invertible matrix with entries
in $R$ and write $\phi^* \e = \bgamma \d \bomega + \bet \bomega$ where $\gamma,
\eta$ are $g \times g$ matrices with entries in $R$. By (7) we have
$\bar{\gamma}= - \bl$.  We get
\[ \begin{array}{lll}
\d^2 \bomega + \balpha \d \bomega + \bbeta \bomega &  = & \psi^*(\bxi \e) =
\bxi \psi^* \e
=\bxi \d \phi^* \e\\
\  & = &  \bxi [ \bgamma \d^2 \bomega + (\d \bgamma + \bet) \d \bomega + (\d
\bet) \bomega ]
\end{array} \]
We get that $\bxi \bgamma = id, \balpha = \bxi \d \bgamma + \bxi \bet, \bbeta =
\bxi \d \bet$. So we have the equations:
\[ \bar{\gamma} \balpha = \d \bar{\gamma} + \bar{\eta},\ \ \ \bar{\gamma}
\bar{\b} = \d \bar{\eta}
\]
Deriving the first equation and  eliminating $\d \bar{\eta}$ we get the
equation (12) and we are done.

\bs

\bs

{\bf Acknowledgements} The first author would like to thank the Institute
for Advanced Study for the hospitality and support during the preparation of
this
paper. He is also grateful to the University of Texas at Austin for an
invitation in June 1994, when the present paper was written.
The authors would also like
to thank the NSF for financial support through grants  DMS 9304580 (A. Buium)
and DMS-9301157 (J. F. Voloch)

\bs

\centerline{\bf References}

\bs

[B1] A.Buium, Intersections in jet spaces and a conjecture of S.Lang, Annals of
Math.
136 (1992), 557-567.

[B2] A.Buium, Geometry of differential polynomial functions,  Part I:   Amer.
J. Math. 115, 6 (1993), 1385-1444; Part II: Amer. J. Math. 116 (1994), 1-33;
Part III: Amer J. Math, to appear.

[B3] A.Buium, Differential Algebraic Groups of Finite Dimension, Lecture Notes
in Math. 1506, Springer 1992.

[BV] A.Buium, J.F.Voloch, The Mordell conjecture in characteristic $p$: an
explicit bound, preprint IAS, 1994.

[Ch] Ching-Li Chai, A note on Manin's Theorem of the Kernel,
Amer. J. Math., 113 (1991), 387-389.

[Co] R.Coleman,  Manin's proof of the Mordell conjecture
over function fields, L'Ens. Math., 36 (1990), 393-427.

[H] E.Hrushovski, The Mordell-Lang conjecture for function fields, preprint,
1993.

[K] N.Katz, Algebraic solutions of differential equations (p-curvature
and the Hodge filtration), Inventiones Math.  18, (1972) 1-118.

[Man1]  Yu.I.Manin, Rational points on algebraic curves
over function fields, Izvestija Akad Nauk SSSR, Mat.Ser.t.27
(1963),1395-1440.

[Man2]  Yu.I.Manin, Algebraic curves over fields
with differentiation, Izv. Akad. Nauk SSSR, Ser. Mat.
22 (1958), 737-756 = AMS translations Series 2,
37 (1964), 59-78.

[Man3] Yu.I.Manin, The Hasse-Witt matrix of an algebraic curve,
AMS translations Series 2, 45 (1965), 245-264.

[MM]  B.Mazur, W.Messing, Universal Extensions and One Dimensional
Crystalline Cohomology, Lecture Notes in Math.370, Springer 1974.

[V1] J.F.Voloch, Explicit $p-$descent for elliptic curves in characteristic
$p$, Compositio Math. 74 (1990), 247-258.

[V2] J.F.Voloch, Diophantine Approximation on Abelian varieties in
in characteristic $p$, Preprint, 1994.

\bs
\bs

\noindent Institute for Advanced Study, Princeton, NJ 08540\\
\ \\
 Dept. of Mathematics, Univ. of Texas,
Austin, TX 78712, USA

\end{document}